\documentclass[12pt]{article}

\usepackage{amssymb,amsmath,amsthm}
\usepackage{amsfonts}
\usepackage{enumerate}
\usepackage{etoolbox}
\usepackage{mathtools}
\usepackage{fancyhdr}	
\usepackage{mdframed}
\usepackage{geometry}
\newtheoremstyle{mytheoremstyle} % name
    {10pt}                    % Space above
    {10pt}                    % Space below
    {\normalfont}                   % Body font
    {}                           % Indent amount
    {\bfseries}                   % Theorem head font
    {.}                          % Punctuation after theorem head
    {0.3cm}                       % Space after theorem head
    {}  % Theorem head spec (can be left empty, meaning ‘normal’)
\theoremstyle{mytheoremstyle}

\theoremstyle{plain}% Theorem-like structures provided by amsthm.sty

\fancypagestyle{plain}{%
  \fancyhf{}%
}

\newtheorem{theorem}{Theorem}[section]

\newtheorem{proposition}[theorem]{Proposition}

\newtheorem{remark}[theorem]{Remark}

\newtheorem{example}[theorem]{Example}

\begin{document}

\title{Quasi-orthogonal extension of skew-symmetric matrices}

\author{Abderrahim Boussa\"{\i}ri, Brahim Chergui, Zaineb Sarir, Mohamed Zouagui}  

\maketitle

\begin{abstract}
 	A real matrix $Q$ is quasi-orthogonal if $Q^{\top}Q=qI$, for some positive real number $q$. We prove that any $n\times n$ skew-symmetric matrix $S$ is a principal sub-matrix of a skew-symmetric quasi-orthogonal matrix $Q$, called a quasi-orthogonal extension of $S$. Moreover, we determine the least integer $d$ such that $S$ has a quasi-orthogonal extension of order $n+d$. This integer is called the quasi-orthogonality index of $S$. Lastly, we give a spectral characterization of skew-adjacency matrices of tournaments with quasi-orthogonality index at most three.
\end{abstract}



\section{Introduction}
Unless specified otherwise, all the matrices considered here are real. We denote by $I_{n}$ and $O_{n}$ respectively the $n \times n$ identity matrix and all-zeros matrix.

An $n\times n$ real matrix $Q$ is \emph{quasi-orthogonal} if $Q^{\top}Q=qI_{n}$, for some positive real number $q$ or equivalently, the matrix $\frac{1}{\sqrt{q}}Q$ is orthogonal. Quasi-orthogonal matrices were first highlighted by Sylvester \cite{sylvester1867lx} and later by Hadamard \cite{hadamard1893resolution}. Well-known classes of quasi-orthogonal matrices are Hadamard matrices and conference matrices. Recall that a \emph{Hadamard matrix} is a square matrix with entries in $\{-1,1\}$ and whose columns are mutually orthogonal \cite{paley1933orthogonal,hedayat1978hadamard,koukouvinos2008skew}. The order of such matrices must be 1, 2 or a multiple of 4. A \emph{conference matrix} is an $n\times  n$ matrix $C$ with 0 on the diagonal and $\pm 1$ elsewhere such that $C^{\top}C=(n-1)I_{n}$ (see for example \cite{belevitch1950theorem,goethals1967orthogonal,belevitch1968conference}). Clearly, an $n\times n$ skew-symmetric matrix $C$ is a  conference matrix if and only if $C+I_{n}$ is a Hadamard matrix. It follows that the order of a skew-symmetric conference matrix is either 2 or multiple of 4. It is conjectured that skew-symmetric conference matrices exist for every order divisible by four \cite{wallis1971some}.

An $n$-\emph{tournament} $T$ is a digraph of order $n$ in which every pair of vertices is joined by exactly one arc. If the arc joining two given vertices $u$ and $v$ of $T$ is directed from $u$ to $v$, we say that $u$ \emph{dominates} $v$. A tournament is \emph{doubly regular} if every pair of distinct vertices dominate the same number $d\geq1$ of vertices.
Brown and Reid \cite{reid1972doubly} proved that the existence of a skew-symmetric conference matrix of order $n = 4k + 4$ is equivalent to the existence of a doubly regular tournament of order $n = 4k + 3$.

Let $T$ be a tournament with vertex set $\{v_{1},\dots,v_{n}\}$, the \emph{skew-adjacency matrix} of $T$ is the skew-symmetric matrix $S_{T}$ in which for $i\neq j\in \{1,\dots,n\}$ the $(i,j)$-entry is $1$ if $v_{i}$ dominates $v_{j}$ and $-1$ otherwise.
Let $T$ be a doubly regular tournament and let $T^{+}$ be the tournament obtained from $T$ by adding a vertex that dominates all vertices of $T$. It follows from the proof of Theorem 2 in \cite{reid1972doubly} that the skew-adjacency matrix of $T^{+}$ is a skew-symmetric conference matrix. Examples of doubly regular tournaments are Paley tournaments. For a prime power $q$ congruent to $3\pmod 4$, the \emph{Paley tournament} of order $q$ is the tournament whose vertex set is the Galois field  $GF(q)$, such that a vertex $x$ dominates a vertex $y$ if $y-x$ is a square in $GF(q)$. Graham and Spencer \cite{graham1971constructive} proved that every tournament can be embedded in a Paley tournament. Therefore, any skew-symmetric matrix with $\pm 1$ off the diagonal is a principal sub-matrix of a skew-symmetric conference matrix.		

In our first main result, we prove that any skew-symmetric matrix $S$ is a principal sub-matrix of a skew-symmetric quasi-orthogonal matrix $Q$. We then say that $Q$ is a \emph{quasi-orthogonal extension} of $S$. The \emph{quasi-orthogonality index} of a skew-symmetric matrix $S$ is the least integer $d$, denoted by ${\rm ind}(S)$, such that $S$ has a quasi-orthogonal extension of order $n+d$. For example,  we will see in Section \ref{index}, that ${\rm ind} (O_{n})=n$.

Let $S$ be a non-zero skew-symmetric matrix of order $n$. The non-zero eigenvalues of $S$ are pure imaginary and come in conjugate pairs.
Then the spectrum of $S$ can be written as $\{[\pm i\lambda_{1}]^{m_{1}},\dots,[\pm i\lambda_{r}]^{m_{r}},[0]^{m}\}$, where $\lambda_{1}>\lambda_{2}>\dots>\lambda_{r}>0$ and $m$ could equal zero. The spectral radius $\rho(S)$ of $S$ is $\lambda_{1}$ and the multiplicity $m_{1}$ of the eigenvalue $i\rho(S)$ will be denoted by $\mu_{S}$ or simply $\mu$ if no confusion occurs. Under these notations, we prove that the quasi-orthogonality index of $S$ is given by the following formula
\[
{\rm ind}(S)=n-2\mu \mbox{.}
\]

In the last section, we investigate on the quasi-orthogonality index of skew-adjacency matrices of tournaments. We define the quasi-orthogonality index of a tournament as the quasi-orthogonality index of its skew-adjacency matrix. We give a spectral characterization of tournaments with quasi-orthogonal index less than four.

A forthcoming article \cite{symmetricMatrices} will be devoted to the quasi-orthogonal extension of symmetric matrices.

\section{Existence of quasi-orthogonal extensions}
For the zero matrix $O_{n}$, the existence of a quasi-orthogonal extension is assured by the given quasi-orthogonal matrix
\[
\begin{pmatrix}
	O_{n}    & I_{n}\\
	-I_{n}  & O_{n}
\end{pmatrix}. 
\]

The existence of quasi-orthogonal extensions for non-zero skew-symmetric matrices is guaranteed by the following theorem.
\begin{theorem}
	\label{existence}
	Let $S$ be a non-zero $n\times n$ skew-symmetric matrix with spectral radius $\rho$. Then $S$ has a quasi-orthogonal extension $\hat{S}$ of order $2n-2\mu$, where $\mu$ is the multiplicity of the eigenvalue $i\rho$.
\end{theorem}
The proof of this theorem  requires the following eigenvalue perturbation result due to Brauer \cite[Theorem 27]{brauer1952limits}.
\begin{theorem}
	\label{Brauer}
	Let $A$ be an $n\times n$ complex matrix with eigenvalues $\lambda,\lambda_{2},\dots,\lambda_{n}$
	and let $\mathbf{x}$ be a non-zero eigenvector of $A$ associated with $\lambda$. Then for any vector $\mathbf{v} \in\mathbb{C}^{n}$, the eigenvalues
	of $A + \mathbf{xv}^{*}$ are $\lambda + \mathbf{v}^{*} \mathbf{x},\lambda_{2},\dots,\lambda_{n}$.
\end{theorem}

\begin{proof}[Proof of Theorem \ref{existence}]
	
	Let $m:=n-2\mu$. 
	
	If $m=0$ then the spectrum of $S$ is $[\pm i\rho]^{\mu}$ and hence the minimal polynomial of $S$ is $x^{2}+\rho^{2}$. Consequently, $S^{2}+\rho^{2}I_{n}=O_n$, which implies that $S$ is quasi-orthogonal.

	If $m=1$ then the spectrum of $S$ is $\{[\pm i\rho]^{\mu}, [0]^{1}\}$. Let
	\[
	\hat{S}:=
	\begin{pmatrix}
		S & \mathbf{u}\\
		-\mathbf{u}^{\top} & 0
	\end{pmatrix}\mbox{,}		
	\]
	where  $\mathbf{u}$ is a vector in the nullspace of $S$ with norm  $\rho$. Then, we have
	
	\[
	\hat{S}^{2}=\begin{pmatrix}
		S^{2}-\mathbf{u}\mathbf{u}^{\top} & S\mathbf{u}\\
		-\mathbf{u}^{\top}S & -\mathbf{u}^{\top}\mathbf{u}
	\end{pmatrix}		
	=
	\begin{pmatrix}
		S^{2}-\mathbf{u}\mathbf{u}^{\top} & \mathbf{0}_{n}\\
		\mathbf{0}_{n}^{\top} & -\rho^{2}
	\end{pmatrix}	\mbox{,}	
	\]
	where $\mathbf{0}_{n}$ is the $n$-dimensional zero vector.
	Since the spectrum of  $S^{2}$ is $\{[-\rho^{2}]^{2\mu}, [0]^{1}\}$, Brauer’s Theorem implies that the spectrum of $S^{2}-\mathbf{u}\mathbf{u}^{\top}$ is  $\{[-\rho^{2}]^{2\mu+1}\}$. Therefore, the spectrum of $\hat{S}^{2}$ is $\{[-\rho^{2}]^{2\mu+2}\}$. It follows that $\hat{S}$ is a quasi-orthogonal extension of $S$ of order $2n-2\mu$. 
	
	We continue the proof by induction on $m$ for $m\geq 2$.
	Let $\alpha$ be the maximum of $|\lambda|$ where $\lambda$ is an eigenvalue of $S$ such that $|\lambda|<\rho$. The real number $\alpha$ is well-defined because $m\neq 0$. We distinguish two cases.  
	
	$\bullet$ \textbf{Case 1:} $\alpha=0$. The spectrum of $S$ is $\{[\pm i\rho]^{\mu}, [0]^{m}\}$.  Then the spectrum of $S^{2}$ is $\{[-\rho^{2}]^{2\mu}, [0]^{m}\}$.  Let $\mathbf{x}$ and $\mathbf{y}$ be two orthogonal vectors in the nullspace of $S^{2}$ with norm $\rho$ and let
	\[
	\hat{S}:=\begin{pmatrix}
		S & \mathbf{x} & \mathbf{y}\\
		-\mathbf{x}^{\top} & 0 & 0\\
		-\mathbf{y}^{\top} & 0 & 0
	\end{pmatrix}.
	\]
	We have
	\begin{equation}\label{brauer2}
		\hat{S}^{2}=\begin{pmatrix}
			S^{2}-\mathbf{x}\mathbf{x}^{\top} -\mathbf{y}\mathbf{y}^{\top} & \mathbf{0}_{n} & \mathbf{0}_{n}\\
			\mathbf{0}_{n}^{\top} & -\mathbf{x}^{\top}\mathbf{x} & 0\\
			\mathbf{0}_{n}^{\top} & 0 & -\mathbf{y}^{\top}\mathbf{y}
		\end{pmatrix}
		\mbox{.}
	\end{equation} 
	Brauer’s theorem applied to the eigenpair $(0,\mathbf{x})$ of $S^{2}$ implies that the spectrum of $S^{2}-\mathbf{xx}^{\top}$ is 
	$\{[-\rho^{2}]^{2\mu +1}, [0]^{m -1}\}$. Since $\mathbf{y}$ is in the  nullspace of $S$ and $\mathbf{x}\bot \mathbf{y}$, we have $(S^{2}-\mathbf{x}\mathbf{x}^{\top})\mathbf{y}=0$. 	
	Reapplying Brauer's theorem to the eigenpair $(0,\mathbf{y})$ of $S^{2}-\mathbf{x}\mathbf{x}^{\top}$, the spectrum of $S^{2}-\mathbf{x}\mathbf{x}^{\top} -\mathbf{y}\mathbf{y}^{\top}$ is $\{[-\rho^{2}]^{2\mu +2}, [0]^{m-2}\}$.
	Using $(\ref{brauer2})$, the spectrum of $\hat{S}^{2}$ is $\{[-\rho^{2}]^{2\mu +4}, [0]^{m-2}\}$ and hence the spectrum of $\hat{S}$ is $\{[\pm i\rho]^{\mu+2},[0]^{m-2}\}$.
	
	$\bullet$ \textbf{Case 2:} $\alpha\neq 0$. The spectrum of  $S$ is  $\{[\pm i\rho]^{\mu}, [\pm i\alpha]^{\tau},\dots \}$ where $\tau \ge 1$. Then the spectrum of $S^{2}$ is $\{[-\rho^{2}]^{2\mu},[-\alpha^{2}]^{2\tau},\dots\}$. Let $\mathbf{x}$ be an eigenvector of $S^{2}$ corresponding to the eigenvalue $-\alpha^{2}$ such that $||\mathbf{x}||^{2}=\rho^{2}-\alpha^{2}$. Let
	\begin{equation*}
		\hat{S}:=\begin{pmatrix}
			S & \mathbf{x} & \frac{1}{\alpha}S\mathbf{x}\\
			-\mathbf{x}^{\top} & 0 & \alpha\\
			\frac{1}{\alpha}\mathbf{x}^{\top}S &-\alpha & 0
		\end{pmatrix} \mbox{.}
	\end{equation*}
	Then
	\begin{equation*}
		\hat{S}^{2}=\begin{pmatrix}
			S^{2}-\mathbf{x}\mathbf{x}^{\top} -\frac{1}{\alpha^{2}}S\mathbf{x}(S\mathbf{x})^{\top} & \mathbf{0}_{n} & \mathbf{0}_{n}\\
			\mathbf{0}_{n}^{\top} & -\mathbf{x}^{\top}x-\alpha^{2} & 0\\
			\mathbf{0}_{n}^{\top} & 0 & -\frac{1}{\alpha^{2}}(S\mathbf{x})^{\top}S\mathbf{x}-\alpha^{2}
		\end{pmatrix} \mbox{.}
	\end{equation*}
	We have
	\begin{align*}
		\frac{1}{\alpha^{2}}(S\mathbf{x})^{\top}S\mathbf{x}&=\frac{1}{\alpha^{2}}\mathbf{x}^{\top}S^{\top}S\mathbf{x}\\
		&=\frac{-1}{\alpha^{2}}\mathbf{x}^{\top}S^{2}\mathbf{x}\\
		&=\frac{\alpha^{2}}{\alpha^{2}}\mathbf{x}^{\top}\mathbf{x}\\
		&=||\mathbf{x}||^{2}=\rho^{2}-\alpha^{2},
	\end{align*}
	and hence
	\begin{equation*}
		\hat{S}^{2}=\left(\begin{array}{ccc}
			S^{2}-\mathbf{x}\mathbf{x}^{\top} -\frac{1}{\alpha^{2}}S\mathbf{x}(S\mathbf{x})^{\top} & \mathbf{0}_{n} & \mathbf{0}_{n}\\
			\mathbf{0}_{n}^{\top} & -\rho^{2} & 0\\
			\mathbf{0}_{n}^{\top} & 0 & -\rho^{2}
		\end{array}
		\right).
	\end{equation*}
	As $(-\alpha^{2},\mathbf{x})$ is an eigenpair of $S^{2}$, Brauer’s theorem implies that the spectrum of $S^{2}-\mathbf{x}\mathbf{x}^{\top}$ is 
	\[\{[-\rho^{2}]^{2\mu},[-\alpha^{2}-\mathbf{x}^{\top}\mathbf{x}]^{1}, [-\alpha^{2}]^{2\tau-1},\dots\}=\{[-\rho^{2}]^{2\mu+1}, [-\alpha^{2}]^{2\tau-1},\dots\}\mbox{.}\] 
	It is easy to check that $(S^{2}-\mathbf{x}\mathbf{x}^{\top})S\mathbf{x}=-\alpha^{2}S\mathbf{x}$, then by Brauer’s theorem, the spectrum of $S^{2}-\mathbf{x}\mathbf{x}^{\top}-\frac{1}{\alpha^{2}}S\mathbf{x}(S\mathbf{x})^{\top}$ is \[\{[-\rho^{2}]^{2\mu+1},[-\alpha^{2}-\frac{1}{\alpha^{2}}(S\mathbf{x})^{\top}S\mathbf{x}]^{1},[-\alpha^{2}]^{2\tau-2},\dots\}= \{[-\rho^{2}]^{2\mu+2},[-\alpha^{2}]^{2\tau-2},\dots\}.\] 
	Therefore, the spectrum of $\hat{S}^{2}$ is $\{[-\rho^{2}]^{2\mu+4}, [-\alpha^{2}]^{2\tau-2},\dots\}$ and hence the spectrum of $\hat{S}$ is $\{[\pm i\rho]^{\mu+2}, [\pm i\alpha]^{\tau-1},\dots\}$.
	
	For $m\geq 2$, we have constructed a quasi-orthogonal extension $\hat{S}$ of $S$ of order $\hat{n}:=n+2$, with spectral radius $\rho$. The multiplicity of $i\rho$ in $\hat{S}$ is $\hat{\mu}:=\mu+2$, then
	\begin{align*}
		\hat{m}&:=\hat{n}-2\hat{\mu}\\
		&=n+2-2(\mu+2)\\
		&= m-2<m
	\end{align*}
	By induction hypothesis, $\hat{S}$ has a quasi-orthogonal extension, which is also a quasi-orthogonal extension of $S$, of order
	\begin{align*}
		2\hat{n}-2\hat{\mu}&=2(n+2)-2(\mu+2)\\
		&=2n-2\mu.
	\end{align*} 	
	
\end{proof}
The strategy used in the proof gives an effective algorithm for embedding a skew-symmetric matrix into a skew-symmetric quasi-orthogonal matrix. In fact, Brauer's theorem produces a bordered matrix in which a conjugate pair
of eigenvalues is increased to the spectral radius. By applying this method repeatedly, we construct a matrix in which all eigenvalues equal the spectral radius. Such a matrix is easily seen to be skew-symmetric and quasi-orthogonal.

%\begin{remark}\label{important_point}
%	Here, we are going to list two important points:
%	\begin{itemize}
	%		\item 	
	%		Notice that the quasi-orthogonal extension of $S$, obtained in the proof of Theorem \ref{existence}, has the same spectral radius as $S$. This is not always the case. For example, consider the matrix 
	%		
	%		\[
	%		S=
	%		\begin{pmatrix*}[r]
		%			0 & 1 & 1 & 1 \\ 
		%			-1 & 0 & 1 & 1  \\ 
		%			-1 & -1 & 0 & -1 \\ 
		%			-1 & -1 & 1 & 0 \\ 
		%		\end{pmatrix*}.
	%		\]
	%		The spectrum of $S$ is $\{[\pm i(\sqrt{2}+1)],[\pm i(\sqrt{2}-1)]\}$.
	%		The matrix 
	%		\[
	%		Q=
	%		\begin{pmatrix*}[r]
		%			${\rm\textbf{0}}$ & ${\rm\textbf{1}}$ & ${\rm\textbf{1}}$ & ${\rm\textbf{1}}$ & -1 & 1 & 1 & 1 \\ 
		%			-${\rm\textbf{1}}$ & ${\rm\textbf{0}}$ & ${\rm\textbf{1}}$ & ${\rm\textbf{1}}$ & -1 & -1 & -1 & -1 \\ 
		%			-${\rm\textbf{1}}$ & -${\rm\textbf{1}}$ & ${\rm\textbf{0}}$ & -${\rm\textbf{1}}$ & -1 & -1 & 1 & 1 \\ 
		%			-${\rm\textbf{1}}$ & -${\rm\textbf{1}}$ & ${\rm\textbf{1}}$ & ${\rm\textbf{0}}$ & 1 & 1 & 1 & -1 \\ 
		%			1 & 1 & 1 & -1 & 0 & -1 & 1 & -1 \\ 
		%			-1 & 1 & 1 & -1 & 1 & 0 & -1 & 1 \\ 
		%			-1 & 1 & -1 & -1 & -1 & 1 & 0 & -1 \\ 
		%			-1 & 1 & -1 & 1 & 1 & -1 & 1 & 0%
		%		\end{pmatrix*}
	%		\]
	%		is a quasi-orthogonal extension of $S$ with spectrum  $\{[\pm i(\sqrt{7})]^{4}\}$.
	%		\item
	%		Quasi-orthogonal matrices are invertible, they have an even order. Hence the quasi-orthogonality index of a matrix has the same parity as its order.
	%	\end{itemize}
%\end{remark} 
\section{Quasi-orthogonality index }\label{index}

As seen in the previous section, $O_n$ has a quasi-orthogonal extension of order $2n$. Then, the quasi-orthogonality index of $O_n$ is at most $n$. 
We now consider a quasi-orthogonal extension $Q$ of $O_n$.
Without loss of generality, we can assume that 
\[Q=
\begin{pmatrix*}[r]
	O_{n}    & B\\
	-B^{\top}  & A
\end{pmatrix*}\mbox{,} 
\]
where $A$ is a square matrix of order $m$ and $B$ is an $n\times m$ matrix. We have 
${\rm rank}\begin{pmatrix}
	O_{n}    & B\\
\end{pmatrix} \leq m$.
Moreover, as $Q$ is invertible, 
${\rm rank} \begin{pmatrix}
	O_{n}    & B\\
\end{pmatrix}=n$. It follows that $m\geq n$ and hence ${\rm ind} (O_{n})\geq n$. 
We conclude  that ${\rm ind} (O_{n})=n$.

The following theorem provides the exact value of the quasi-orthogonality index for non-zero skew-symmetric matrices.
\begin{theorem}	\label{exact_value}
	Let $S$ be an $n\times n$ non-zero skew-symmetric matrix with spectral radius $\rho$. Then the quasi-orthogonality index of $S$ is $n-2\mu$ where $\mu$ is the multiplicity of the eigenvalue $i\rho$. Moreover, if $\hat{S}$ is a  quasi-orthogonal extension of $S$ of order
	$2n-2\mu$ then the spectral radius of $\hat{S}$ is $\rho$. 
\end{theorem}
To prove this theorem, we need the following proposition, which is a direct consequence of
Cauchy's Interlace Theorem \cite{cauchyInterlacing}.

\begin{proposition}\label{sameeigenvalue}
	Let $H$ be an $n\times n$ Hermitian matrix and let $H^{\prime}$ be a principal submatrix of $H$ of order $m$. Then any eigenvalue of $H$ with multiplicity $r$ more than $n-m$ is  an eigenvalue of $H^{\prime}$ with multiplicity at least $r-n+m$.
\end{proposition}
\begin{proof}[Proof of Theorem \ref{exact_value}]
	By Theorem \ref{existence}, $S$ has a quasi-orthogonal extension of order $2n-2\mu$. Then ${\rm ind}(S)\leq n-2\mu$. Let $\hat{S}$ be a quasi-orthogonal extension of $S$ of order $\hat{n}:=n+{\rm ind}(S)$. By Remark \ref{important_point}, we have $\hat{n}=2n-2\mu-2h$ for some non negative integer $h$.
	We denote by $\hat{\rho}$ the spectral radius of $\hat{S}$. Consider the Hermitian matrices $H=iS$ and $\hat{H}=i\hat{S}$.
	Clearly, $\hat{\rho}$ is an eigenvalue of $\hat{H}$ with multiplicity $\mu_{\hat{S}}=n-\mu-h$. In addition, $\mu_{\hat{S}}>\hat{n}-n$. By Proposition \ref{sameeigenvalue}, $\hat{\rho}$ is an eigenvalue of $H$ with multiplicity at least $\mu +h$. Then  $\hat{\rho}=\rho$ and hence $h=0$.
	
\end{proof}
%\label{majorindex}
\begin{remark}\label{important_point}
	Here, we are going to list two useful points:
	\begin{itemize}
		\item Skew-symmetric quasi-orthogonal matrices are invertible, so they have an even order. Hence, the quasi-orthogonality index of a matrix has the same parity as its order.
		\item 	By Theorem \ref{exact_value}, the quasi-orthogonality index of a non-zero $n\times n$ skew-symmetric matrix is at most $n-2$.
	\end{itemize}
\end{remark}

\begin{example}
	Consider the following matrix
	\[
	S=
	\begin{pmatrix*}[r]
		0 & 1 & 1 & 1 \\ 
		-1 & 0 & 1 & 1  \\ 
		-1 & -1 & 0 & -1 \\ 
		-1 & -1 & 1 & 0 \\ 
	\end{pmatrix*}.
	\]
	The spectrum of $S$ is $\{[\pm i(\sqrt{2}+1)],[\pm i(\sqrt{2}-1)]\}$. Then, by Theorem \ref{exact_value}, $ind(S)=2$.
	The technique used in the proof of Theorem \ref{existence} allows us to construct the following quasi-orthogonal extension of $S$
	$$\begin{pmatrix*}[r]
		0 & 1 & 1 & 1 & 2^{3/4} & 0 \\
		-1 & 0 & 1 & 1 & -\sqrt[4]{2} & -\sqrt[4]{2} \\
		-1 & -1 & 0 & -1 & \sqrt[4]{2} & -\sqrt[4]{2} \\
		-1 & -1 & 1 & 0 & 0 & 2^{3/4} \\
		-2^{3/4} & \sqrt[4]{2} & -\sqrt[4]{2} & 0 & 0 & \sqrt{2}-1 \\
		0 & \sqrt[4]{2} & \sqrt[4]{2} & -2^{3/4} & 1-\sqrt{2} & 0 \\
	\end{pmatrix*}.$$
\end{example}

\section{Quasi-orthogonality index of tournaments}
In this section, we give a spectral characterization of tournaments with small quasi-orthogonality index. First, notice that 
the quasi-orthogonality index is 0 for the unique 2-tournament and 1 for the two 3-tournaments.

Let $T$ be an $n$-tournament with skew-adjacency matrix $S_{T}$. We define the skew-characteristic polynomial $\phi_{T}(x)$ of $T$ as the characteristic polynomial of $S_{T}$.

The following theorem provides a characterization of tournaments with quasi-orthogonality index less than four in terms of their characteristic polynomials.
\begin{theorem}\label{Tournament_index}
	Let $T$ be an $n$-tournament with $n\geq 4$. Then, we have the following statements
	\begin{enumerate}[i)]
		\item If $n=4k+4$, then ${\rm ind}(T)\geq0$. Equality holds if and only if $\phi_{T}(x)=(x^{2}+4k+3)^{2k+2}$.
		\item If $n=4k+3$, then ${\rm ind}(T)\geq1$. Equality holds if and only if $\phi_{T}(x)=x(x^{2}+4k+3)^{2k+1}.$
		\item If $n=4k+2$, then ${\rm ind}(T)\geq2$. Equality holds if and only if $\phi_{T}(x)=(x^{2}+4k+3)^{2k}(x^{2}+1)$.
		\item If $n=4k+1$, then ${\rm ind}(T)\geq3$. Equality holds if and only if either $n=5$ or $n\geq 9$ and $\phi_{T}(x)=x(x^{2}+4k+3)^{2k -1}(x^{2}+3)$.
	\end{enumerate}
\end{theorem}
Before proving Theorem \ref{Tournament_index}, we recall some results about tournaments.  
Let $\phi_{T}(x)=x^{n}+a_{1}x^{n-1}+\cdots+a_{n}$ be the characteristic polynomial of an $n$-tournament $T$. Then 
\begin{equation}\label{coefficients}
	a_{k}=(-1)^{k}\sum(\text{all } k\times k\text{ principal minors}). 
\end{equation}
In particular, we have
\begin{align}
	a_{2}&=\binom{n}{2},\label{minors}\\ 
	a_{n}&=\det(S_{T}),\\
	a_{k}&=0 \ \text{if} \ k \ \text{is odd}.
\end{align}

\begin{proposition}[\cite{McCarthy1996}]\label{determinant}
	If $n$ is even then the determinant of $S_{T}$ is an odd integer.
\end{proposition}
\begin{proposition}[\cite{McCarthy1996}]\label{nullSpace}
	The nullspace of $S_{T}$ has dimension zero if $n$ is even, and dimension one if $n$ is odd.
\end{proposition}
The proof of Theorem \ref{Tournament_index} also requires some algebraic properties of eigenvalues of tournaments.

Let $\alpha$ be an algebraic number with minimal polynomial $m_{\alpha}(x)$. The norm of $\alpha$, denoted by ${\rm N}(\alpha)$, is the product of all the roots of $m_{\alpha}(x)$, which is equal  up to sign to the constant coefficient of $m_{\alpha}(x)$. In particular, if $\alpha$ is an algebraic integer then ${\rm N}(\alpha)$ is an integer.

\begin{proposition}\label{oddNorm}
	Let $S_{T}$ be a skew-adjacency matrix of an $n$-tournament $T$, and $i\lambda$ be a non-zero eigenvalue of $S_{T}$. Then $i\lambda$ is an algebraic integer and ${\rm N}(i\lambda)$ is odd.
\end{proposition}
\begin{proof}
	As $S_{T}$ is an integral matrix, $\phi_{T}(x)$ is a monic polynomial with integral coefficients and hence $i\lambda$ is an algebraic integer.
	If $n$ is even, then ${\rm N}(i\lambda)$ divides $\det(S_{T})$ because $m_{i\lambda}(x)$ divides $\phi_{T}(x)$. By Proposition \ref{determinant}, $\det(S_{T})$ is odd, thus ${\rm N} (i\lambda)$ is odd. If $n$ is odd, then $m_{i\lambda}(x)$ divides the polynomial $\phi_{T}(x)/x$. It follows that ${\rm N} (i\lambda)$ divides the constant coefficient of $\phi_{T}(x)/x$, which is, by \eqref{coefficients}, equal to the sum of the $n$ principal minors of order $n-1$. Due to Proposition \ref{determinant}, this sum is odd. Consequently ${\rm N} (i\lambda)$ is odd. 
\end{proof}

\begin{proposition}\label{integercoefficient}
	Let $P(x)=(x^{2}+\rho^{2})^{l}(x^{2}+\alpha^{2})$ where $\alpha$, $\rho$ are real numbers such that $\alpha\neq\pm \rho$ and $l>1$. If $P(x)$ has integral coefficients, then $\alpha^{2}$ and $\rho^{2}$ are integers.
\end{proposition}
\begin{proof}
	The complex numbers $i\rho$ and $i\alpha$ are algebraic integers. Moreover, the minimal polynomial of $i\rho$ is either $(x^{2}+\rho^{2})(x^{2}+\alpha^{2})$ or $x^{2}+\rho^{2}$.
	Suppose that $m_{i\rho}(x)=(x^{2}+\rho^{2})(x^{2}+\alpha^{2})$. As $m_{i\rho}(x)$ is monic polynomial with integral coefficients, which divides $P(x)$, then $P(x)/m_{i\rho}(x)= (x^{2}+\rho^{2})^{l-1}$ is a monic polynomial with integral coefficients that annihilates $i\rho$, contradiction.
	It follows that $m_{i\rho}(x)=x^{2}+\rho^{2}$, in particular $\rho^{2}$ is an integer. The polynomial $(x^{2}+\rho^{2})^{l}$ is monic with integral coefficients, then $P(x)/(x^{2}+\rho^{2})^{l}=x^{2}+\alpha^{2}$ has integral coefficients and hence $\alpha^{2}$ is an integer.
	
\end{proof}

\begin{proof}[Proof of Theorem \ref{Tournament_index}]
	Let $\mu$ be the multiplicity of $i\rho$, where $\rho$ is the spectral radius of $S_{T}$.
	We have 
	\begin{equation}
		1\leq \mu\leq \lfloor n/2\rfloor 
	\end{equation}
	\begin{enumerate}[i)]
		
		\item If ${\rm ind}(T)=0$, then $S_T$ is a quasi-orthogonal matrix. By definition, there exists a positive real number $q$ such that $S_{T}^2+qI_n=O_n$. The minimal polynomial $\pi_{S_T}$ of $S_T$ is $\pi_{S_T}(x)=x^2+q$. The eigenvalues of $S_{T}$ are $\pm i\sqrt{q}$. Then $q=\rho^2$ and $\phi_{T}(x)=(x^2+\rho^{2})^{2k+2}$. By \eqref{minors}, we get $(2k+2)\rho^{2}=(4k+4)(4k+3)/2$ and then $\rho^{2}=4k+3$. It follows that $\phi_{T}(x)=(x^{2}+4k+3)^{2k+2}$. 
		
		Conversely, If $\phi_{T}(x)=(x^{2}+4k+3)^{2k+2}$ then $\pi_{S_{T}}(x)=x^2+4k+3$, which implies that $S_T$ is a quasi-orthogonal matrix and ${\rm ind}(T)=0$. 
		
		\item As $n$ is odd, by Remark \ref{important_point}, ${\rm ind}(T)$ is odd and hence ${\rm ind}(T)\geq1$.
		
		If ${\rm ind}(T)=1$ then, by Theorem $\ref{exact_value}$, $\mu=2k+1$. Moreover, as $S_{T}$ is skew-symmetric of odd order, $0$ is an eigenvalue of $S_{T}$. Thus $\phi_{T}(x)=x(x^{2}+\rho^{2})^{2k+1}$. The same argument used in the proof of i) gives $\rho^{2}=4k+3$. Hence, $\phi_{T}(x)=x(x^{2}+4k+3)^{2k+1}$.
		
		Conversely, if $\phi_{T}(x)=x(x^{2}+4k+3)^{2k+1}$ then, by Theorem \ref{exact_value}, ${\rm ind}(T)=1$.
		
		\item If ${\rm ind}(T)=0$ then $S_{T}$ is a quasi-orthogonal matrix and $\phi_{T}(x)=(x^{2}+\rho^{2})^{2k+1}$. By \eqref{minors}, we have $\rho^{2}=n-1$. Then, the minimal polynomial of $S_{T}$ is $x^{2}+n-1$ and hence  $S_{T}^{2}+(n-1)I_n=O_{n}$. It follows that $S_{T}$ is a skew-conference matrix, which contradicts the fact that the order of a skew-conference matrix is 2 or divisible by $4$. Then, by Remark \ref{important_point}, ${\rm ind}(T)\geq 2$.
		
		If ${\rm ind}(T)=2$ then, by Theorem \ref{exact_value}, $\mu=2k$. Since $n$ is even, by Proposition \ref{nullSpace}, $\phi_{T}(x)=(x^{2}+\rho^{2})^{2k}(x^{2}+\alpha^{2})$, where $\alpha$ is a positive real number not equal to $\rho$.
		By \eqref{minors}, we have
		\begin{equation}\label{minor'}
			4k\rho^{2} + 2\alpha^{2}=(4k+2)(4k+1). 
		\end{equation}	
		It follows that 
		\begin{equation*}
			\rho^{2} \leq 4k+3+\dfrac{1}{2k}\mbox{.}
		\end{equation*}
		Moreover, by \cite[Corollary 2.3.]{Chen2015}, the skew spectral radius of an $n$-tournament is at least equal to $\sqrt{n-1}$. Then
		\[4k+1\leq \rho^{2}\leq 4k+3 +\dfrac{1}{2k}.\]
		According to Proposition \ref{integercoefficient}, $\rho^{2}$ and $\alpha^{2}$ are integers. Then, ${\rm N}(i\rho)=\rho^{2}$ and ${\rm N}(i\alpha)=\alpha^{2}$. By Proposition \ref{oddNorm}, $\rho^{2}$ and $\alpha^{2}$ are odd. It follows that $\rho^{2}\in\{4k+1,4k+3\}$. If $\rho^{2}=4k+1$ then by \eqref{minor'} we get $\alpha^{2}=4k+1$, which is a contradiction because $\alpha\neq\rho$. We then must have $\rho^{2}=4k+3$. Therefore, $\alpha^{2}=1$. Hence $\phi_{T}(x)=(x^{2}+4k+3)^{2k}(x^{2}+1)$. The converse follows from Theorem \ref{exact_value}.
		
		\item By Remark \ref{important_point}, ${\rm ind}(T)$ is odd. If ${\rm ind}(T)=1$ then, by Theorem \ref{exact_value}, $\mu=2k$. As seen in the proof of ii), $\phi_{T}(x)=x(x^{2}+\rho^{2})^{2k}$. By \eqref{minors}, we have $\rho^{2}=4k+1$. Let $\mathbf{u}:=\begin{pmatrix}
			u_{1}\\
			u_{2}\\
			\vdots\\
			u_{n}
		\end{pmatrix}$ be an eigenvector of $S_{T}$ associated with the eigenvalue $0$ such that $\lVert  \mathbf{u}\lVert=\rho$. 
		According to the proof of Theorem \ref{existence}, the matrix
		\[
		\hat{S}:=\begin{pmatrix}
			S_{T} & \mathbf{u}\\
			-\mathbf{u}^{\top} & 0
		\end{pmatrix}
		\]
		is a quasi-orthogonal extension of $S_{T}$. More precisely, $\hat{S}^{2}+\rho^{2}I_{4k+2}=O_{4k+2}$. The norm of each column vector in $\hat{S}$ is equal to $4k + u_{i}^{2}=\rho^{2}=4k+1$. Hence $u_{i}=\pm 1$. Then $\hat{S}$ is a skew-conference matrix of order $4k+2$. This contradicts the fact that the order of a skew-conference matrix is 2 or divisible by 4. It follows from Remark \ref{important_point} that ${\rm ind}(T)\geq3$.

		Suppose that ${\rm ind}(T)=3$ and $n\geq 9$. By Theorem \ref{exact_value}, we have $\mu=2k-1$. It follows from Proposition \ref{nullSpace} that $\phi_{T}(x)=x(x^{2}+\rho^{2})^{2k -1}(x^{2}+\alpha^{2})$, where $\alpha$ is a positive real number not equal to $\rho$. 
		The polynomial $\phi_{T}(x)$ and, a fortiori, $\phi_{T}(x)/x$ have integral coefficients. Since $2k-1=\dfrac{n-3}{2}>1$, Proposition \ref{integercoefficient} implies that $\alpha^{2}$ and $\rho^{2}$ are integers and hence $\alpha^{2}\geq 1$.
		By \eqref{minors} we have 
		\begin{align}\label{rhoAlpha}
			(2k-1)\rho^{2} + \alpha^{2}=4k(4k+1)/2.
		\end{align}
		Then
		\begin{align*}
			\rho^{2}=\frac{2k(4k+1)-\alpha^{2}}{2k-1}.
		\end{align*}		
		It follows that 
		\begin{align*}
			\rho^{2} \leq \frac{2k(4k+1)-1}{2k-1}.
		\end{align*}
		As 
		\[\frac{2k(4k+1)-1}{2k-1}=4k+3+\frac{2}{2k-1},\]
		we get $\rho^{2}\leq 4k+3$, because $k\geq 2$ and $\rho^{2}$ is an integer.
		By using \cite[Corollary 2.3.]{Chen2015} we have $4k\leq \rho^{2}\leq 4k+3$. Since $\rho^{2}$ and $\alpha^{2}$ are integers, ${\rm N}(i\rho)=\rho^{2}$ and ${\rm N}(i\alpha)=\alpha^{2}$. By Proposition \ref{oddNorm}, $\rho^{2}$ and $\alpha^{2}$ are odd and hence $\rho^{2}\in\{4k+1,4k+3\}$.
		If $\rho^{2}=4k+1$, then, by \eqref{rhoAlpha}, $\alpha^{2}=4k+1$, which is a contradiction because $\alpha\neq\rho$.
		Then $\rho^{2}=4k+3$, which leads to $\alpha^{2}=3$ and hence $\phi_{T}(x)=x(x^{2}+4k+3)^{2k -1}(x^{2}+3)$.

		If $n\geq 9$ and $\phi_{T}(x)=x(x^{2}+4k+3)^{2k -1}(x^{2}+3)$ then by Theorem \ref{exact_value}, we have ${\rm ind}(T)=3$. 
		
		If $n=5$ then by Remark \ref{important_point}, we have also ${\rm ind}(T)=3$.
		
	\end{enumerate}
\end{proof}

It follows from assertion i) of Theorem \ref{Tournament_index} that the quasi-orthogonality index of a $n$-tournament is zero if and only if its skew-adjacency matrix is a skew-symmetric conference matrix. As mentioned in the introduction, the existence of skew-symmetric conference of order $n$ is an open problem. If such a matrix exists then $n$ divisible by four, moreover the following theorem ensures the existence of tournaments with order $n-1$, $n-2$ and $n-3$ whose quasi-orthogonality indices are respectively 1, 2 and 3.
\begin{theorem}[\cite{greaves2017symmetric}]
	The existence of the following are equivalent: 
	\begin{enumerate}[i)]
		\item a skew-symmetric Seidel \footnote{Following \cite{greaves2017symmetric}, a \emph{Seidel matrix} is a $\{0, \pm1\}$-matrix $S$ with zero diagonal and all off-diagonal entries nonzero such that $S = \pm S^{\top}$.} matrix with characteristic polynomial $$(x^{2} + 4k + 3)^{2k+2};$$ \item a skew-symmetric Seidel matrix with characteristic polynomial $$x(x^{2} + 4k + 3)^{2k+1};$$ \item a skew-symmetric Seidel matrix with characteristic polynomial $$(x^{2} + 1)(x^{2} + 4k + 3)^{2k};$$
		\item  a skew-symmetric Seidel matrix with characteristic polynomial $$x(x^{2} + 3)(x^{2} + 4k + 3)^{2k-1}.$$
	\end{enumerate}
\end{theorem}
\bibliography{quasibib}

\end{document}